\documentstyle[12pt,titlepage]{article}

\begin{document}
\begin{titlepage}
\vspace*{\fill}
\begin{center}
\large\bf ASYMPTOTIC APPROXIMATIONS FOR SYMMETRIC ELLIPTIC
INTEGRALS
\end{center}
\vfill
\begin{center}
\large B. C. Carlson\(\dagger\) and John L. Gustafson\(\ddagger\)
\end{center}
\vfill
\begin{center}
Dedicated to Dick Askey and Frank Olver in gratitude for many years
of
friendship.
\end{center}
\vfill
\hyphenation{Le-gen-dre}
{\bf Abstract.}\quad Symmetric elliptic integrals, which have been
used as
replacements for Legendre's integrals in recent integral tables and
computer
codes, are homogeneous functions of three or four variables.  When
some of the
variables are much larger than the others, asymptotic
approximations with
error bounds are presented.  In most cases they are derived from a
uniform
approximation to the integrand.  As an application the symmetric
elliptic
integrals of the first, second, and third kinds are proved to be
linearly
independent with respect to coefficients that are rational
functions.
\vfill
{\bf Key words.}\quad Elliptic integral, asymptotic approximation,
inequalities, hypergeometric \(R\)-function
\\
\\
{\bf AMS(MOS) subject classifications.}\quad primary 33A25, 41A60
26D15;
secondary 33A30, 26D20
\\
\\
This work was supported by the Director of Energy Research, Office
of Basic
Energy Sciences.  The Ames Laboratory is operated for the
U.~S.~Department of
Energy by Iowa State University under Contract W-7405-ENG-82.

\(\dagger\)\, Ames Laboratory and Department of Mathematics, Iowa
State
 University, Ames, Iowa 50011-3020.

\(\ddagger\)\, Ames Laboratory and Department of Computer Science,
Iowa State
 University, Ames, Iowa 50011-3020.
\end{titlepage}
\vspace*{\fill}
\begin{center}
Abbreviated \vspace{1.0in} title\\
{\sc\bf ASYMPTOTIC APPROXIMATIONS FOR ELLIPTIC INTEGRALS}
\end{center}
\vfill
\newpage
\section{Introduction}
A recent table of elliptic integrals \cite{C87,C88,C89,C91,C92}
uses symmetric standard integrals instead of Legendre's
integrals because permutation symmetry makes it possible to unify
many of the
formulas in previous tables.  Fortran codes for numerical
computation of the
symmetric integrals, which are homogeneous functions of three or
four variables,
can be found in several major software libraries as well as in the
supplements
to \cite{C87,C88}. For analytical purposes it is desirable to
know how the homogeneous functions behave when some of the
variables are much
larger than the others.  For all such cases we list in Section 2
asymptotic
approximations (sometimes two or three approximations of different
accuracy),
always with error bounds.  Proofs are discussed in Section 3.  In
most cases
the approximations are obtained by replacing the integrand by a
uniform
approximation.  Many of the results found by a different method
in \cite{G} have been improved by sharpening the error
bounds or by finding bounds for incomplete elliptic integrals that
are still
useful for the complete integrals, which are then not listed
separately.  Cases
not considered in \cite{G} include two for a completely symmetric
integral of
the second kind and two for a symmetric integral of the third kind
in which
two variables are much larger than the other two.

We assume that \(x,y,z\) are nonnegative and at most one of them is
0.  The
symmetric integral of the first kind,
\begin{equation}\label{1.1}
R_F(x,y,z) =
\frac{1}{2}\int_{0}^{\infty}[(t+x)(t+y)(t+z)]^{-1/2}dt,
\end{equation}
is homogeneous of degree \(-1/2\) in \(x,y,z\) and satisfies
\(R_F(x,x,x) =
x^{-1/2}\).
The symmetric integral of the third kind,
\begin{equation}\label{1.2}
R_J(x,y,z,p) =
\frac{3}{2}\int_{0}^{\infty}[(t+x)(t+y)(t+z)]^{-1/2}(t+p)^{-1}dt,
\quad p>0,
\end{equation}
is homogeneous of degree \(-3/2\) in \(x,y,z,p\) and satisfies
\(R_J(x,x,x,x) =
x^{-3/2}\).  If \(p=z,\, R_J\) reduces to an integral of the second
kind,
\begin{equation}\label{1.3}
R_D(x,y,z) = R_J(x,y,z,z) =
\frac{3}{2}\int_{0}^{\infty}[(t+x)(t+y)]^{-1/2}(t+z)^{-3/2}dt,
\quad z>0,
\end{equation}
which is symmetric in \(x\) and \(y\) only.  If two variables of
\(R_F\)
 are equal, the integral becomes an elementary function,
\begin{equation}\label{1.4}
R_C(x,y) = R_F(x,y,y) =
\frac{1}{2}\int_{0}^{\infty}(t+x)^{-1/2}(t+y)^{-1}dt, \quad y>0.
\end{equation}
If \(x<y\) it is an inverse trigonometric function,
\begin{equation}\label{1.5}
R_C(x,y) = (y-x)^{-1/2}\arccos(x/y)^{1/2},
\end{equation}
and if \(x>y\) it is an inverse hyperbolic function,
\begin{equation}\label{1.6}
R_C(x,y) = (x-y)^{-1/2}{\rm arccosh}(x/y)^{1/2} =
(x-y)^{-1/2}\ln\frac{\sqrt{x}+\sqrt{x-y}}{\sqrt{y}}.
\end{equation}
If the second argument of \(R_C\) is negative, the Cauchy principal
value is
\cite[(4.8)]{Z}
\begin{equation}\label{1.7}
R_C(x,-y) = \left(\frac{x}{x+y}\right)^{1/2} R_C(x+y,y), \quad y>0.
\end{equation}
If the fourth argument of \(R_J\) is negative, the Cauchy principal
value is
given by \cite[(4.6)]{Z}
\begin{eqnarray}\
(y+p)R_J(x,y,z,-p) & = & (q-y)R_J(x,y,z,q) - 3R_F(x,y,z)
\nonumber\\
& + & 3\left(\frac{xyz}{xz+pq}\right)^{1/2}R_C(xz+pq,pq), \quad
p>0,
\label{1.8}
\end{eqnarray}
where \(q-y = (z-y)(y-x)/(y+p)\).  If we permute the values of
\(x,y,z\) so
 that \( (z-y)(y-x) \geq 0 \), then \(q \geq y > 0 \).

A completely symmetric integral of the second kind is not as
convenient as
\(R_D\) for use in tables because its representation by a single
integral is
more complicated \cite[(9.1-9)]{C77}:
\begin{equation}\label{1.9}
R_G(x,y,z) = \frac{1}{4}\int_{0}^{\infty}[(t+x)(t+y)(t+z)]^{-1/2}
\left( \frac{x}{t+x} + \frac{y}{t+y} + \frac{z}{t+z} \right) tdt.
\end{equation}
It is symmetric and homogeneous of degree \(1/2\) in \(x,y,z\), and
it
satisfies \(R_G(x,x,x) = x^{1/2}\).  It has a
nice representation by a double integral that expresses the surface
area of an
ellipsoid \cite[(9.4-6)]{C77}.  It is related to \(R_D\)
and \(R_F\) by (\ref{3.4}) and by
\begin{eqnarray}
6R_G(x,y,z)\!&=&\!x(y+z)R_D(y,z,x) + y(z+x)R_D(z,x,y)
+ z(x+y)R_D(x,y,z),\quad \label{1.11} \\
6R_G(x,y,0) &=& xy[R_D(0,x,y) + R_D(0,y,x)].\hfill \label{1.12}
\end{eqnarray}

Legendre's complete elliptic integrals \(K\) and \(E\) are given by
\begin{eqnarray}
K(k) & = & R_F(0,1-k^2,1), \label{1.13}\\
E(k) & = & 2R_G(0,1-k^2,1) \nonumber\\
& = & \frac{1-k^2}{3} [R_D(0,1-k^2,1) + R_D(0,1,1-k^2)],
\label{1.14}\\
K(k) - E(k) & = & \frac{k^2}{3} R_D(0,1-k^2,1), \label{1.15}\\
E(k) - (1-k^2)K(k)  & = & \frac{k^2(1-k^2)}{3} R_D(0,1,1-k^2).
\label{1.16}
\end{eqnarray}
Approximations and
inequalities for $K$, $E$, and some combinations thereof are given
in
\cite{V0,V,B}.  If the error terms in (\ref{F1e}),\,(\ref{F1f}),
and
(\ref{G1c}) are omitted, the approximations reduce to the leading
terms of
well-known series expansions of $K$ and $E$ for $k$ near $1$
\cite[p.~54]{C}
\cite[900.06, 900.10]{BF}.  If the series for $K$ is truncated
after any number
of terms, simple bounds for the {\em relative} error are given in
\cite[(1.17)]{CG}.  A generalization of this series to
\(R_F(x,y,z)\) with
\(x,y <\!< z\) is given in \cite[(1.14)-(1.16)]{CG}, again with
simple bounds
for the relative error of truncation.

The various functions designated by \(R\) with a letter subscript
are special
cases of the multivariate hypergeometric \(R\)-function,
\[
R_{-a}(b_1,\ldots,b_n;\,z_1,\ldots,z_n),
\]
which is symmetric in the indices \(1,\ldots,n\) and homogeneous of
degree
\(-a\) in the variables \(z_1,\ldots,z_n\).  Best regarded as the
Dirichlet
average of \(x^{-a}\) \cite[\S~5.9]{C77}, it is a symmetric variant
of the
function known as Lauricella's \(F_D\).  By the method of Mellin
transforms,
series expansions are obtained in \cite[(4.16)-(4.19)]{C85} that
converge
rapidly if some of the \(z\)'s are much larger than the others and
if the
parameters satisfy \( \sum_{i=1}^{n}b_i > a > 0 \).  Thus
the leading terms of these series provide asymptotic approximations
for all
except \(R_G\) among the functions
\def \hf{\textstyle{\frac{1}{2}}}
\def \th{\textstyle{\frac{3}{2}}}
\begin{eqnarray}
R_F(x,y,z) & = & R_{-\frac{1}{2}}(\hf,\hf,\hf;\,x,y,z), \quad 
R_C(x,y) =
R_{-\frac{1}{2}}(\hf,1;\,x,y), \label{1.17}\\
R_J(x,y,z,p) & = & R_{-\frac{3}{2}}(\hf,\hf,\hf,1;\,x,y,z,p), \quad
R_D(x,y,z) = R_{-\frac{3}{2}}(\hf,\hf,\th;\,x,y,z),\label{1.18}\\
R_G(x,y,z) & = & R_{\frac{1}{2}}(\hf,\hf,\hf;\,x,y,z). \label{1.19}
\end{eqnarray}
However, error bounds for the approximations are more easily
derived by the
methods of the present paper.  Another function that is used
repeatedly in
obtaining error bounds is \cite[Ex.~9.8-5]{C77}
\begin{eqnarray}
R_{-1}(\hf,\hf,1;\,x,y,z) & = &
\int_{0}^{\infty}[(t+x)(t+y)]^{-1/2}
 (t+z)^{-1}dt \label{1.20}\\
 & = & 2R_C((\sqrt{xy}+z)^2,(\sqrt{x}+\sqrt{y})^2z).\label{1.21}
\end{eqnarray}

In Section 4 the asymptotic approximations are applied to show that
\(R_F(x,y,z)\), \(R_D(x,y,z), R_J(x,y,z,p)\), and \( (xyz)^{-1/2}
\)
are linearly independent with respect to coefficients that
are rational functions of
\(x,y,z\), and \(p\).  An Appendix contains some elementary
inequalities that are used in obtaining error bounds.

The results in this paper provide upper and lower approximations
that approach
the elliptic integrals as selected ratios of the variables approach
zero.
Approximations that approach the integrals as all variables
approach a common
value have been found by other methods.  For example, the theory of
hypergeometric mean
values yields upper and lower algebraic approximations for all the
integrals in
this paper \cite[Thm. 2]{C66}, while truncation of Taylor series
about the
arithmetic mean of the variables gives approximations with errors
that may be
positive or negative.  Successive applications of the duplication
theorem for
\(R_F\), making its three variables approach equality, provide
ascending and
descending sequences of successively sharper (and successively more
complicated) algebraic approximations to \(R_F\) and \(R_C\)
\cite{C70}.
Transcendental approximations that approach \(R_F\) when only two
of its
variables approach equality are furnished by
\begin{equation}\label{1.22}
   R_C\!\left( x,\frac{y+z}{2} \right) \leq R_F(x,y,z)
   \leq R_C(x,\sqrt{yz}),\qquad yz \neq 0,
\end{equation}
which follows from (\ref{AX}).  The inequalities can be sharpened
by first
using Landen or Gauss transformations of \(R_F\) \cite[\S~9.5]{C77}
to make
$y$ and $z$ approach equality.  If \(x = 0\) the Gauss
transformation reduces
to replacing \( \sqrt{y} \) and \( \sqrt{z} \) by their arithmetic
and geometric means, and each \(R_C\)-function
becomes \( \pi /2 \) divided by the square root of its second
argument.  Therefore, in the complete case the procedure reduces to
the
algorithm of the arithmetic-geometric mean
\cite[(6.10-6)(9.2-3)]{C77}
and provides ascending and
descending sequences of algebraic approximations, of which leading
members are shown in (\ref{F2b}).
\section{Results}
We assume throughout that \(x,y\), and \(z\) are nonnegative and at
most one of
them is 0.  The last argument of \(R_C\), \(R_D\), and \(R_J\) is
assumed to be
positive (see (\ref{1.7}) and (\ref{1.8})).
\\
\\
\noindent {\bf C1.}{\boldmath \(\quad R_C(x,y)\ \)}{\bf with}
{\boldmath \ \(x<\!<y.\)}
\begin{equation}\label{C1}
R_C(x,y) = \frac {\pi}{2\sqrt{y}} - \frac{\sqrt{x}}{y}
 + \frac{\pi x \theta}{4y^{3/2}}\,,
\end{equation}
where \(1/(1+\sqrt{x/y}) \leq \theta \leq 1\) with equalities iff
\(x=0\).
\\
\\
\noindent {\bf C2.}{\boldmath \(\quad R_C(x,y)\ \)}{\bf with}
{\boldmath \(\ y<\!<x.\ \)} Two approximations of different
accuracy are
\begin{eqnarray}
R_C(x,y) & = & \frac{1}{2\sqrt{x}} \left(\ln\frac{4x}{y} +
\frac{y}{2x-y}
\ln\frac{\theta_1 x}{y} \right) \label{C2a}\\
& = & \frac{1}{2\sqrt{x}} \left[
\left(1+\frac{y}{2x}\right)\ln\frac{4x}{y}
- \frac{y}{2x} + \frac{3y^2}{4x(2x-y)}\ln\frac{\theta_2 x}{y}
\right],
\label{C2b}
\end{eqnarray}
where \(1< \theta_i < 4\) for $i=1,2$.  The first approximation
implies
\begin{equation}\label{C2c}
R_C(x,y) < \frac{1}{2\sqrt{x}(1-y/2x)} \ln\frac{4x}{y}.
\end{equation}
\\
\\
\noindent {\bf F1.}{\boldmath \(\quad R_F(x,y,z)\ \)}{\bf with}
{\boldmath \(\ x,y<\!<z.\ \)} Let \(a = (x+y)/2\,,\, g =
\sqrt{xy}\,\), and
\(\rho = max\{x,y\}/z\).  Then
\begin{equation}\label{F1a}
R_F(x,y,z) = \frac{1}{2\sqrt{z}}\left(\ln\frac{8z}{a+g} +
\frac{r}{2z}\right),
\end{equation}
where
\[
\frac{g}{1-g/z}\ln\frac{2z}{a+g} < r <
\frac{a}{1-a/2z}\ln\frac{8z}{a+g}.
\]
The upper bound implies
\begin{equation}\label{F1b}
R_F(x,y,z) < \frac{1}{2\sqrt{z}(1-a/2z)}\ln\frac{8z}{a+g}.
\end{equation}
A sharper lower bound and a higher-order approximation are given by
\begin{eqnarray}
R_F(x,y,z) & = & \frac{1}{2\sqrt{z}} \left(\ln\frac{8z}{a+g}
 + \frac{ar_1}{2z} \right)\label{F1c}\\
& = & \frac{1}{2\sqrt{z}}
\left[\left(1+\frac{a}{2z}\right)\ln\frac{8z}{a+g}
 - \frac{2a-g}{2z} + \frac{3(3a^2 - g^2)r_2}{16z^2}
\right],\label{F1d}
\end{eqnarray}
where
\[
\ln\frac{z}{2a} < \frac{\ln(1/\rho)}{1-\rho} < r_i <
\frac{1}{1-a/2z}\ln\frac{8z}{a+g}\,,\quad i = 1,2.
\]
By (\ref{1.13}) this implies (since $4k^2 < 4-k'^2$ if $k^2 < 1$)
\begin{eqnarray}
K(k) & = & \ln\frac{4}{k'} +
\frac{k'^2}{4-k'^2}\ln\frac{\theta_1}{k'}\label{F1e}\\
& = & \left(1+\frac{k'^2}{4}\right)\ln\frac{4}{k'} - \frac{k'^2}{4}
 + \frac{9k'^4}{16(4-k'^2)}\ln\frac{\theta_2}{k'}, \label{F1f}
\end{eqnarray}
where \(0 < k' = \sqrt{1-k^2}\) and \(1 < \theta_i < 4\) for \(i =
1,2\).
\\
\\
\noindent {\bf F2.}{\boldmath \(\quad R_F(x,y,z)\ \)}{\bf with}
{\boldmath \(\ z<\!< x,y.\ \)} Let \(a = (x+y)/2\) and \(g =
\sqrt{xy}\,\). Then
\begin{equation}\label{F2a}
R_F(x,y,z) = R_F(x,y,0) - \frac{\sqrt{z}}{g} + \frac{\pi
z\theta}{4g^{3/2}},
\end{equation}
where \(1/(1+\sqrt{z/g}) < \theta < a/g\).  Note that
\( R_F(x,y,0)= \pi /2AGM(\sqrt{x},\sqrt{y})\),
where \(AGM\) denotes Gauss's arithmetic-geometric mean
\cite[(6.10-6)(9.2-3)]{C77}, and hence
\begin{equation}\label{F2b}
\frac{1}{\sqrt{a}} \leq
\sqrt{\frac{2}{a+g}} \leq \frac{2}{\sqrt{(a+g)/2}+\sqrt{g}} \leq
\frac{2}{\pi}R_F(x,y,0) \leq \left(\frac{2}{ag+g^2}\right)^{1/4}
\leq
\frac{1}{\sqrt{g}}\,,
\end{equation}
with equalities iff $x = y$.
\\
\\
\noindent {\bf D1.}{\boldmath \(\quad R_D(x,y,z)\ \)}{\bf with}
{\boldmath \(\ x,y<\!<z.\ \)} Let \(a = (x+y)/2\) and \( g =
\sqrt{xy}\,\). Then
\begin{equation}\label{D1}
R_D(x,y,z) = \frac{3}{2z^{3/2}}\left(\ln\frac{8z}{a+g} - 2
 + \frac{\theta}{z}\ln\frac{2z}{a+g} \right),
\end{equation}
where
\[
\frac{g}{1-g/z} < \theta < \frac{3a}{2(1-a/z)}\,.
\]
\\
\\
\noindent {\bf D2.}{\boldmath \(\quad R_D(x,y,z)\ \)}{\bf with}
{\boldmath \(\ z<\!<x,y.\ \)} Let \(a = (x+y)/2\) and \( g =
\sqrt{xy}\,\). Then
\begin{equation}\label{D2a}
R_D(x,y,z) = \frac{3}{\sqrt{xyz}}\left(1 -
\frac{\pi\theta}{2}\sqrt{\frac{z}{g}}
 \right),
\end{equation}
where
\[
1 - \frac{4}{\pi}\sqrt{\frac{z}{g}} < \theta < \frac{a}{g}.
\]
A higher-order approximation is
\begin{equation}\label{D2b}
R_D(x,y,z) = \frac{3}{\sqrt{xyz}} - R_D(0,x,y) - R_D(0,y,x)
 + \frac{3\pi\theta \sqrt{z}}{2g^2(1+\sqrt{z/g})},
\end{equation}
where
\[
\frac{1}{\sqrt{2/3}+\sqrt{z/g}} < \theta <
\frac{3a}{2g(1+\sqrt{z/g})}.
\]
An approximation of still higher order is
\begin{equation}\label{D2c}
R_D(x,y,z) = \frac{3}{\sqrt{xyz}} - \frac{6}{xy}R_G(x,y,0)
 + \frac{6a\sqrt{z}}{g^3}\left(1-\frac{\pi
\theta}{4}\sqrt{\frac{z}{a}}\right),
\end{equation}
where we have used (\ref{1.12}) and where
\[
\frac{1}{1+\sqrt{z/a}} < \theta < \left(\frac{a}{g}\right)^{3/2}
\left(3-\frac{g^2}{a^2} \right),
\]
\\
\\
\noindent {\bf D3.}{\boldmath \(\quad R_D(x,y,z)\ \)}{\bf with}
{\boldmath \(\ y,z<\!<x.\ \)} Let \(a = (y+z)/2\) and \( g =
\sqrt{yz}\,\). Then
\begin{equation}\label{D3}
R_D(x,y,z) = \frac{3}{\sqrt{x}}\left(\frac{1}{g+z} - \frac{r}{4x}
\right),
\end{equation}
where
\[
\frac{1}{1-g/x}\ln\frac{2x}{a+g} - \frac{2z}{g+z} < r <
\frac{1}{1-a/2x}
\ln\frac{8x}{a+g}.
\]
\\
\\
\noindent {\bf D4.}{\boldmath \(\quad R_D(x,y,z)\ \)}{\bf with}
{\boldmath \(\ x<\!<y,z.\ \)} Let \(a = (y+z)/2\) and \( g =
\sqrt{yz}\,\). Then
\begin{equation}\label{D4}
R_D(x,y,z) = R_D(0,y,z) + \frac{3\sqrt{x}}{gz}\left(-1 + \frac{\pi
\theta}{4}
 \sqrt{\frac{x}{a}} \right),
\end{equation}
where
\[
\frac{1}{1+\sqrt{x/a}} < \theta < \left(\frac{a}{g} \right)^{3/2}
 \left(1+\frac{y}{a} \right).
\]
\\
\\
\noindent {\bf J1.}{\boldmath \(\quad R_J(x,y,z,p)\ \)}{\bf with}
{\boldmath \(\ x,y,z<\!<p.\ \)} Let \(a = (x+y+z)/3\) and
 \( b = (\sqrt{3}/2)(xy+xz+yz)^{1/2}\,\). Then
\begin{equation}\label{J1a}
R_J(x,y,z,p) = \frac{3}{p}R_F(x,y,z) + \frac{3\pi}{2p^{3/2}}(-1+r),
\end{equation}
where
\[
\frac{\sqrt{b/p}}{1+\sqrt{b/p}} < r <
\frac{3}{2}\ \frac{\sqrt{a/p}}{1+\sqrt{a/p}}.
\]
In the complete case a sharper result is
\begin{equation}\label{J1b}
R_J(x,y,0,p) = \frac{3}{p}\left(R_F(x,y,0) -
\frac{\pi}{2\sqrt{p}}\right)
 \left(1 + \frac{\theta/p}{1-\theta /p}\right),
\end{equation}
where \(\sqrt{xy} \leq \theta \leq (x+y)/2\) with equalities iff
$x=y$.
\\
\\
\noindent {\bf J2.}{\boldmath \(\quad R_J(x,y,z,p)\ \)}{\bf with}
{\boldmath \(\ p<\!<x,y,z.\ \)} Let \(g = (xyz)^{1/3}\),
 \(3h^{-1} = x^{-1}+y^{-1}+z^{-1}\), and \(\lambda =
 \sqrt{xy}+\sqrt{xz}+\sqrt{yz}\). Note that \(g\) is the geometric
mean and
\(h\) is
the harmonic mean, whence $g \geq h$ with equality iff $x=y=z$.
Then
\begin{equation}\label{J2a}
R_J(x,y,z,p) = \frac{3}{2\sqrt{xyz}}\left(\ln\frac{4g}{p} -2 +
r\right),
\end{equation}
where
\[
-\ln\frac{g}{h} < r < \frac{3p}{2(g-p)}\ln\frac{g}{p}.
\]
A higher-order approximation is
\begin{equation}\label{J2b}
R_J(x,y,z,p) = \frac{3}{2\sqrt{xyz}}\ln\frac{4xyz}{p\lambda^2}
 + 2R_J(x+\lambda,y+\lambda,z+\lambda,\lambda) +
\frac{3pr}{4\sqrt{xyz}},
\end{equation}
where
\[
\frac{2}{g-p}\ln\frac{g}{p} < r < \frac{3}{h-p}\ln\frac{h}{p}.
\]
The second term in the approximation is independent of \(p\) but is
otherwise
as complicated as the function being approximated.  The same is
true of an even
more accurate approximation \cite[Thm.~11]{G} in which the error is
of order
$p$ instead of $p\ln p$ and the leading term involves $R_C$.
\\
\\
\noindent {\bf J3.}{\boldmath \(\quad R_J(x,y,z,p)\ \)}{\bf with}
{\boldmath \(\ x,y<\!<z,p.\ \)} Let \(a = (x+y)/2\) and \(g =
\sqrt{xy}\). Then
\begin{equation}\label{J3}
R_J(x,y,z,p) = \frac{3}{2\sqrt{z}p}\left[\ln\frac{8z}{a+g} -
 2R_C\left(1,\frac{p}{z}\right) +
\frac{\theta}{p}\ln\frac{2p}{a+g}\right],
\end{equation}
where
\[
\frac{g}{1-g/p} < \theta <
\frac{a}{1-a/p}\left(1+\frac{p}{2z}\right).
\]
\\
\\
\noindent {\bf J4.}{\boldmath \(\quad R_J(x,y,z,p)\ \)}{\bf with}
{\boldmath \(\ z,p<\!<x,y.\ \)} Let \(a = (x+y)/2\), \(g =
\sqrt{xy}\),
\(b=\sqrt{3p(p+2z)}/2\), and \(d = (z+2p)/3\). Then
\begin{equation}\label{J4a}
R_J(x,y,z,p) = \frac{3}{g}R_C(z,p)
 - \frac{3\theta}{g-p}\left[R_C(z,g)-\frac{p}{g}R_C(z,p)\right],
\end{equation}
where \(1 \leq \theta \leq a/g\) with equalities iff $x = y$. Since
$z<\!<g$,
$R_C(z,g)$ can be estimated from (\ref{C1}).  In the complete case
(\ref{J4a})
reduces to
\begin{equation}\label{J4b}
R_J(x,y,0,p) = \frac{3\pi}{2\sqrt{xyp}}\left(1-
 \frac{\theta\sqrt{p}}{\sqrt{g}+\sqrt{p}}\right)
\end{equation}
with $\theta$ as before.  A higher-order approximation is
\begin{equation}\label{J4c}
R_J(x,y,z,p) = \frac{3}{g}R_C(z,p) - \frac{6}{xy}R_G(x,y,0)
 + \frac{3\pi \theta}{2xy},
\end{equation}
where we have used (\ref{1.12}) and where
\[
\frac{\sqrt{b}}{1+\sqrt{b/g}} < \theta < \frac{3a}{2g}\>
\frac{\sqrt{d}}{1+\sqrt{d/g}}.
\]
\\
\\
\noindent {\bf J5.}{\boldmath \(\quad R_J(x,y,z,p)\ \)}{\bf with}
{\boldmath \(\ x<\!<y,z,p.\ \)} Let \(a = (y+z)/2\) and \(g =
\sqrt{yz}\). Then
\begin{equation}\label{J5}
R_J(x,y,z,p) = R_J(0,y,z,p) + \frac{3\sqrt{x}}{gp}\left(-1+
 \frac{\pi \theta}{4}\sqrt{\frac{x}{g}}\right),
\end{equation}
where
\[
\frac{\sqrt{g/a}}{1+\sqrt{x/a}} < \theta <
\frac{a}{g}+\frac{g}{p}\,.
\]
\\
\\
\noindent {\bf J6.}{\boldmath \(\quad R_J(x,y,z,p)\ \)}{\bf with}
{\boldmath \(\ y,z,p<\!<x.\ \)} Let \(a = (y+z)/2\) and \(g =
\sqrt{yz}\). Then
\begin{equation}\label{J6a}
R_J(x,y,z,p) = \frac{3}{\sqrt{x}}\left[R_C((g+p)^2,2(a+g)p)
 - \frac{r}{4}\right],
\end{equation}
where
\[
\frac{1}{x-g}\ln\frac{2x}{a+g} - \frac{2p}{x} R_C((g+p)^2,2(a+g)p)
< r
 < \frac{1}{x-a/2}\ln\frac{8x}{a+g}.
\]
In the complete case this reduces to
\begin{equation}
R_J(x,0,z,p) = \frac{3}{\sqrt{xp}}R_C(p,z) - \frac{3s}{4x^{3/2}},
\end{equation}
where
\[
\ln\frac{4x}{z} - 2\sqrt{p}R_C(p,z) < s <
\frac{1}{1-z/4x}\ln\frac{16x}{z}.
\]
\\
\\
\noindent {\bf G1.}{\boldmath \(\quad R_G(x,y,z)\ \)}{\bf with}
{\boldmath \(\ x,y<\!<z.\ \)} Let \(a = (x+y)/2\) and \(g =
\sqrt{xy}\). Then
\begin{equation}\label{G1a}
R_G(x,y,z) = \frac{\sqrt{z}}{2}\left(1+\frac{r}{2z}\right),
\end{equation}
where
\[
\frac{a+g}{2}\ln\frac{2z}{a+g} + 2g - \frac{4a}{3} < r
 < (3a-g)\ln\frac{2z}{a+g} + 2g - \frac{a}{3}.
\]
In the right-hand inequality it is assumed that $5a < z$.  A
sharper result for
the complete case is
\begin{equation}\label{G1b}
R_G(0,y,z) = \frac{\sqrt{z}}{2} + \frac{y}{8\sqrt{z}}
 \left(\ln\frac{16z}{y}-1+ \frac{ys}{2z}\right),
\end{equation}
where
\[
\frac{3}{4}\ln\frac{z}{y} < s <
\frac{1}{1-y/z}\left(\ln\frac{16z}{y} -
\frac{13}{6}\right).
\]
By (\ref{1.14}) this follows from
\begin{equation}\label{G1c}
E(k) = 1 + \frac{k'^2}{2}\left(\ln\frac{4}{k'} - \frac{1}{2} + k'^2
r\right),
\end{equation}
where $0 < k' = \sqrt{1-k^2} <\!< 1$ and
\[
\frac{3}{8}\ln\frac{1}{k'} < r <
\frac{1}{k(1+k)}\left(\ln\frac{4}{k'} -
\frac{13}{12}\right).
\]
\\
\\
\noindent {\bf G2.}{\boldmath \(\quad R_G(x,y,z)\ \)}{\bf with}
{\boldmath \(\ z<\!<x,y.\ \)} Let \(a = (x+y)/2\) and \( g =
\sqrt{xy}\). Then
\begin{equation}\label{G2}
R_G(x,y,z) = R_G(x,y,0) + \pi \theta z/8,
\end{equation}
where
\[
\frac{1}{\sqrt{a}}\left(1 - \frac{4}{\pi}\sqrt{\frac{z}{a}}\right)
< \theta
 < \left(\frac{2}{ag+g^2}\right)^{1/4} \leq \frac{1}{\sqrt{g}}
\]
with equality iff x = y.
\section{Proofs}
Most of the results in Section 2 are obtained by replacing an
integrand
$f$ by an approximation $f_a$, writing
 \(\int f = \int f_a + \int(f-f_a)\),
and finding upper and lower bounds for $\int(f-f_a)$.  All
integrals are taken
over the positive real line.  The function $f_a$ is usually chosen
to be
a uniform approximation
\(f_a = f_i + f_o - f_m\), where $f_i$ is an approximation in the
inner region, $f_o$ in the outer region, and
$f_m$ in the overlap region or matching region.  For
instance, if \(f(t) = [(t+x)(t+y)(t+z)]^{-1/2}\) with $x,y<\!<z$,
we get $f_i$
 by
neglecting $t$ compared to $z$, $f_o$ by neglecting $x$ and $y$
compared to
$t$, and $f_m$ by doing both.  A first example of this process is
the proof
of Lemma 1.
\\
\\
{\bf Lemma 1.}\quad
\begin{em}
If $x \geq 0$, $y \geq 0$, and $0 < x+y <\!< z$, then
\begin{equation}\label{3.1}
\int_{0}^{\infty}\frac{dt}{\sqrt{(t+x)(t+y)}(t+z)} =
\frac{1}{z-\theta}
\ln\frac{2z}{a+g},
\end{equation}
where \(\sqrt{xy} = g \leq \theta \leq a = (x+y)/2\) with
equalities iff
$x=y$.
\end{em}
\\
\\
{\bf Proof.} Let
\begin{eqnarray*}
f(t) & = & \frac{1}{\sqrt{(t+x)(t+y)}(t+z)},\quad f_i(t)
 = \frac{1}{z\sqrt{(t+x)(t+y)}},\\
f_o(t) & = & \frac{1}{t(t+z)},\qquad \qquad f_m(t) = \frac{1}{zt}\,
.
\end{eqnarray*}
Taking \(f_a = f_i+f_o-f_m\), we find
\[
\int_{0}^{\infty}f_a(t)dt = \frac{1}{z}\ln\frac{2z}{a+g}
\]
and
\[
f-f_a = \frac{t}{z(t+z)}\left(\frac{1}{t} -
\frac{1}{\sqrt{(t+x)(t+y)}}\right).
\]
Inequality (\ref{A5}) in the Appendix implies
\[
f-f_a = \frac{\theta}{z\sqrt{(t+x)(t+y)}(t+z)}, \qquad g \leq
\theta \leq a,
\]
and thus
\[
\int f = \int f_a + \int (f-f_a) = \int f_a + \frac{\theta}{z}\int
f
 = \frac{1}{1-\theta /z} \int f_a =
\frac{1}{z-\theta}\ln\frac{2z}{a+g}.\qquad
\Box
\]

As a second example, in which Lemma 1 is used, consider
$R_F(x,y,z)$ with
$x,y<\!<z$.  Let
\begin{eqnarray*}
f(t) & = & \frac{1}{\sqrt{(t+x)(t+y)(t+z)}},
 \quad f_i(t) = \frac{1}{\sqrt{(t+x)(t+y)z}},\\
f_o(t) & = & \frac{1}{t\sqrt{t+z}}, \qquad \qquad
f_m(t) = \frac{1}{\sqrt{z}t}.
\end{eqnarray*}
Taking \(f_a = f_i + f_o - f_m\), we find (with $a$ and $g$ the
same as
before)
\[
\int_{0}^{\infty}f_a(t)dt = \frac{1}{\sqrt{z}}\ln \frac{8z}{a+g}
\]
and
\[
f-f_a = \left(\frac{1}{\sqrt{z}} - \frac{1}{\sqrt{t+z}}\right)
        \left(\frac{1}{t} - \frac{1}{\sqrt{(t+x)(t+y)}}\right).
\]
Inequalities (\ref{A2}) and (\ref{A5}) imply
\[
\frac{g}{2\sqrt{z(t+x)(t+y)}(t+z)} < f-f_a <
 \frac{a}{2z\sqrt{(t+x)(t+y)(t+z)}}.
\]
Hence, by Lemma 1,
\[
\frac{g}{2\sqrt{z}(z-g)}\ln\frac{2z}{a+g} < \int (f-f_a) <
\frac{a}{2z}
 \int f < \frac{a/2z}{1-a/2z} \int f_a \,,
\]
where the last inequality follows from the next to last.  We
complete the proof
of (\ref{F1a}) by noting that
\[
2R_F(x,y,z) = \int f = \int f_a + \int(f-f_a).
\]

Equations (\ref{F1c}) and (\ref{F1d}) are obtained from
\cite[(2.15)(3.25)]{CG} with $w = \infty$.  To derive (\ref{F2a})
we construct
\(f_a = f_i + f_o - f_m\) as usual and find bounds for
$\int(f-f_a)$ by
using(\ref{A1}) and (\ref{A6}).  To simplify the upper bound we
note that
\(R_F(x,y,z) \leq R_F(x,y,0)\) and use (\ref{F2b}).

Equations (\ref{C1}),(\ref{C2a}),(\ref{C2b}), and (\ref{C2c})
follow from
(\ref{F2a}),(\ref{F1a}),(\ref{F1d}), and (\ref{F1b}), respectively,
by
replacing \(x\) by \(y\), replacing \(z\) by \(x\), and
simplifying.

Among the approximations for $R_D$ we need discuss only (\ref{D2a})
and
(\ref{D2c}), since (\ref{D1}), (\ref{D2b}), (\ref{D3}), and
(\ref{D4}) follow
from (\ref{J3}), (\ref{J4c}), (\ref{J6a}), and (\ref{J5}),
respectively, by
putting \(p = z\) and simplifying.  To prove (\ref{D2a}) we let
\[
f(t) = \frac{1}{\sqrt{(t+x)(t+y)}(t+z)^{3/2}}, \qquad f_i(t) =
 \frac{1}{g(t+z)^{3/2}}\,,
\]
choose $f_a = f_i\,$, and apply (\ref{A6}) to get
\begin{eqnarray*}
\frac{t}{g(t+g)(t+z)^{3/2}} & \leq & f_a-f \leq
 \frac{at}{g^2\sqrt{(t+x)(t+y)}(t+z)^{3/2}}\,, \\
\frac{1}{g(g-z)\sqrt{t+z}}\left(\frac{g}{t+g} -
\frac{z}{t+z}\right)
 & \leq & f_a-f < \frac{a}{g^2\sqrt{t+z}(t+g)}\,, \\
\frac{2}{g-z}\left[R_C(z,g) - \frac{\sqrt{z}}{g}\right] & \leq &
\int(f_a-f) <
 \frac{2a}{g^2}R_C(z,g).
\end{eqnarray*}
Use of (\ref{C1}) completes the proof.  Approximation (\ref{D2c})
follows
from applying (\ref{D4}) to two terms on the right side of
\begin{equation}\label{3.2}
R_D(x,y,z) = 3(xyz)^{-1/2} - R_D(z,x,y) - R_D(z,y,x),
\end{equation}
an identity that comes from \cite[(5.9-5)(6.8-15)]{C77}.

In discussing approximations for $R_J$, we define
\[
f(t) = \frac{1}{\sqrt{(t+x)(t+y)(t+z)}(t+p)}
\]
and construct $f_i$, $f_o$, and $f_m$ for each case in the manner
described
at the beginning of this Section. For example, if $x,y,z <\!< p$,
then $f_i$ is
obtained by neglecting $t$ compared to $p$.  Unless otherwise
stated, we define
\(f_a = f_i + f_o - f_m\), take \(\int f_a\) as an approximation to
\(\int f\), and find bounds for \(\int(f-f_a)\) by using the
inequalities in
the Appendix.

To prove (\ref{J1a}) we use (\ref{A9}).  To prove (\ref{J1b}) we
use (\ref{A5})
and note that \(\int(f-f_a) = (\theta /p)\int f\).  Before
discussing
(\ref{J2a}), we consider (\ref{J2b}), in which the error bounds are
easily
found by using (\ref{A10}).  Finding \(\int f_a\) requires an
integration by
parts and a formula of which we omit the proof,
\begin{equation}\label{3.3}
\int_{0}^{\infty}(ln~t)\frac{d}{dt}[(t+x)(t+y)(t+z)]^{-1/2}dt =
 \frac{1}{\sqrt{xyz}}\ln\frac{\lambda^2}{4xyz} -
\frac{4}{3}R_J(x+\lambda,
 y+\lambda,z+\lambda,\lambda),
\end{equation}
where \(\lambda = \sqrt{xy}+ \sqrt{xz}+\sqrt{yz}\).  To have a
simpler
approximation (\ref{J2a}), we define \(f_a = f_i + f_o - f_m\) and
\(\phi_a = f_i + f_s - f_m\), where \(f_o\) has been replaced by
\[
f_s(t) = \frac{1}{t(t+g)^{3/2}}
\]
Then
\[
\int\phi_a = \frac{1}{\sqrt{xyz}}\left(\ln\frac{4g}{p} - 2\right),
\]
and an upper bound for \(\int(f-\phi_a)\) is found by using
 \(\sqrt{(t+x)(t+y)(t+z)} \geq (t+g)^{3/2}\) and (\ref{A4}).  To
find a lower
bound, we note that \(f-f_a > 0\), whence
\[
f-\phi_a = f-f_a + f_o-f_s > f_o-f_s.
\]
A lower bound for \(\int(f_o-f_s)\) follows from (\ref{AZ}).

The straightforward proof of (\ref{J3}) uses (\ref{A5}),
(\ref{A7}), and
Lemma 1.  For the elementary approximation (\ref{J4a}) we choose
\(f_a =f_i\)
and use (\ref{A6a}).  For the more accurate approximation
(\ref{J4c}) we take
\(f_a = f_i + f_o - f_m\) and evaluate \(\int f_a\) by integrating
by parts.
The error bounds follow from (\ref{A6a}) and (\ref{A9}) with two
variables
equated.  To find the error bounds for (\ref{J5}), we use
(\ref{A8}),
(\ref{A1}), and (\ref{AX}) to prove
\[
\frac{\sqrt{t}}{(t+x)(t+a)} < \frac{2gp}{x}(f-f_a) <
 \left(\frac{a}{g} + \frac{g}{p}\right) \frac{1}{\sqrt{t+z}(t+g)}.
\]
After integrating, (\ref{C1}) is used to complete the proof.  In
the case of
(\ref{J6a}), where \(\int(f_o - f_m)\) is infinite, we choose \(f_a
= f_i\)
and evaluate \(\int f_a\) by (\ref{1.21}).  It follows from
(\ref{A2}) that
\[
\frac{1}{2\sqrt{x(t+y)(t+z)}}\left(\frac{1}{t+x} -
\frac{p}{x(t+p)}\right)
 < f_a - f < \frac{1}{2x\sqrt{(t+x)(t+y)(t+z)}},
\]
where we have replaced \(t/(t+p)\) by \(1\) in the upper bound and
\(x/(x-p)\)
by 1 in the lower bound.  We then use (\ref{1.21}), (\ref{3.1}),
and
(\ref{F1b}).

The function $R_G$ can be expressed in terms of $R_F$ and $R_D$ by
 (\ref{1.18}) and \cite[Table 9.3-1]{C77}:
\begin{equation}\label{3.4}
2R_G(x,y,z) = zR_F(x,y,z) - \frac{1}{3}(z-x)(z-y)R_D(x,y,z) +
 \sqrt{\frac{xy}{z}}.
\end{equation}
Applying (\ref{F1a}) and (\ref{D1}), we obtain (\ref{G1a}).  The
error bounds
have been substantially simplified by using the numerical value of
$ln~2$
and assuming \(5a < z\) in the upper bound.  It is not hard to
obtain
 (\ref{G1c}) from a well-known infinite series \cite[p.~54]{C} for
$E(k)$ by
using the inequality
\[
1 + \frac{3}{8}k'^2 < \, _2 F_1(\hf,\th ;2;k'^2) < (1-k'^2)^{-1/2}
= 1/k,
\qquad 0 < k' < 1,
\]
for the hypergeometric function $_2F_1$.  Unfortunately (\ref{3.4})
does not
lead to simple error bounds for (\ref{G2}).  Instead, we define
\(f(z) = R_G(x,y,z)\) and find from \cite[(5.9-9)(6.8-6)]{C77} that
\[
f'(z) =
\frac{1}{8}\int_{0}^{\infty}\frac{tdt}{\sqrt{(t+x)(t+y)}(t+z)^{3/
2}}.
\]
Since this is a strictly decreasing function of $z$, the mean value
theorem
yields \(f(z) = f(0) + zf'(\zeta)\) where
\[
f'(z) < f'(\zeta) < f'(0) = \frac{1}{4}R_F(x,y,0).
\]
By (\ref{AX}) and (\ref{1.5}) we see that
\[
f'(z) \geq \frac{1}{8}
\int_{0}^{\infty}\frac{tdt}{(t+a)(t+z)^{3/2}}
 = \frac{1}{4(a-z)}[-\sqrt{z} + aR_C(z,a)].
\]
Use of (\ref{F2b}) and (\ref{C1}) completes the proof of
(\ref{G2}).
\section{Application to linear independence}
In \cite[Thm. 9.2-1]{C77} it is shown that \(R_F(x,y,z)\),
\(R_G(x,y,z)\), an
integral of the third kind called \(R_H(x,y,z,p)\), and the
algebraic function
\((xyz)^{-1/2}\) are linearly independent with respect to
coefficients that are rational functions of \(x,y,z,p\).
It then follows \cite[\S 9.2]{C77} that
every elliptic integral can be expressed in terms of $R_F$, $R_G$,
$R_H$, and
elementary functions. From (\ref{3.4}) and a known relation
expressing $R_H$ in
terms of $R_J$ and $R_F$, we may conclude that every elliptic
integral can be
expressed in terms of $R_F$, $R_D$, $R_J$, and elementary
functions.  In order
to reach the same conclusion without invoking $R_G$ and $R_H$, we
shall use
the results of this paper to prove the linear independence of
\(R_F, R_D,
R_J\), and \( (xyz)^{-1/2} \) with respect to coefficients that are
rational
functions.
\\
\\
\noindent {\bf Theorem 1.}\quad
\begin{em}
The functions \(R_F(x,y,z)\), \(R_D(x,y,z)\),
\(R_J(x,y,z,p)\), and \( (xyz)^{-1/2} \) are linearly independent
with respect
to coefficients that are rational functions of \(x,y,z\), and $p$.
\end{em}
\\
\\
\noindent {\bf Proof.}\quad Let \(\alpha, \beta, \gamma\), and
\(\delta\) be
rational functions of \(x,y,z\), and \(p\).  We need to prove that
\begin{equation}\label{4.1}
\alpha R_F(x,y,z) + \beta R_D(x,y,z) + \gamma R_J(x,y,z,p) +
\delta\,
(xyz)^{-1/2} \equiv 0
\end{equation}
iff \(\alpha, \beta, \gamma\), and \(\delta\) are identically $0$.
We may assume that these coefficients are polynomials since we can
multiply
all terms by the denominator of any rational function.
As \(p \rightarrow 0\), (\ref{J2a}) shows that \(R_J(x,y,z,p)\)
involves
\( \ln p\) while all other quantities are polynomials in
$p$, whence \(\gamma \equiv 0\).  As \(z \rightarrow \infty\) we
have
\[
   \alpha = az^m(1+O(1/z)),\qquad \beta = bz^n(1+O(1/z)),
\]
where \(m\) and \(n\) are nonnegative integers and \(a\) and \(b\)
are
polynomials in $x,y$, and $p$.  Using (\ref{F1a}) and (\ref{D1})
and
multiplying all terms by \(2z^{3/2}\), we find
\[
   az^{m+1} \left[ \ln \frac{8z}{a+g} + O \left( \frac{\ln z}{z}
\right)
   \right] + 3bz^n \left[ \ln \frac{8z}{a+g} -2 + O \left(
\frac{\ln z}{z}
   \right) \right] + 2\delta\, (xy)^{-1/2}z \equiv 0.
\]
Cancellation of the leading terms in \( \ln z\) requires \(az^{m+1}
+ 3bz^n
\equiv 0\), implying \(n = m+1\) and \(a \equiv  -3b\) and leaving
\[
   O(z^m \ln z) - 6bz^{m+1} + 2\delta\, (xy)^{-1/2}z \equiv 0.
\]
Because the second term is of different order from the first and
does not
have a square root in common with the third, it follows that \(b
\equiv0\),
whence also \(a \equiv 0\).  Since the leading terms of the
polynomials
\(\alpha\) and \(\beta\) are identically $0$, so too are \(\alpha\)
and
\(\beta\).  Finally, with only one term remaining in (\ref{4.1}),
we have
\(\delta \equiv 0\).\qquad \(\Box\)

It is an open question whether Theorem 1 is still true if the
coefficients are
algebraic functions instead of rational functions.  However,
polynomial
coefficients suffice (see the first paragraph of
\cite[\S~9.2]{C77}) to prove
that every elliptic integral can be expressed in terms of \(R_F,
R_D, R_J\),
and elementary functions.
%
\appendix
\begin{center}
\Large \bf Appendix
\end{center}
\section*{Elementary inequalities}
Assuming $x$, $y$, $z$, and $t$ are positive, we list and prove
some
inequalities that are used in this paper to obtain error bounds:
\begin{equation}\label{A1}
\frac{x}{2\sqrt{t}(t+x)} < \frac{1}{\sqrt{t}} -
\frac{1}{\sqrt{t+x}} <
\frac{x}{2t\sqrt{t+x}},
\end{equation}
\begin{equation}\label{A2}
\frac{t}{2\sqrt{x}(t+x)} < \frac{1}{\sqrt{x}} -
\frac{1}{\sqrt{t+x}} <
\frac{t}{2x\sqrt{t+x}},
\end{equation}
\begin{equation}\label{A3}
\frac{1}{t^{3/2}} - \frac{1}{(t+x)^{3/2}} = \frac{\theta
x}{t^{3/2}(t+x)},
\qquad 1 < \theta < \frac{3}{2},
\end{equation}
\begin{equation}\label{A4}
\frac{1}{x^{3/2}} - \frac{1}{(t+x)^{3/2}} = \frac{\theta
t}{x^{3/2}(t+x)},
\qquad 1 < \theta < \frac{3}{2}.
\end{equation}
In the next five inequalities let \(a = (x+y)/2\) and \(g =
\sqrt{xy}\).
Inequalities become equalities in (\ref{A5}), (\ref{A6}), and
(\ref{A6a})
iff $x = y$.
\begin{equation}\label{A5}
\frac{1}{t} - \frac{1}{\sqrt{(t+x)(t+y)}} =
\frac{\theta}{t\sqrt{(t+x)(t+y)}},
\qquad  g \leq \theta \leq a,
\end{equation}
\begin{equation}\label{A6}
\frac{t}{g(t+g)} \leq \frac{1}{\sqrt{xy}} -
\frac{1}{\sqrt{(t+x)(t+y)}}
 \leq \frac{at}{g^2\sqrt{(t+x)(t+y)}},
\end{equation}
or alternatively,
\begin{equation}\label{A6a}
\frac{1}{\sqrt{xy}} - \frac{1}{\sqrt{(t+x)(t+y)}} = \frac{\theta
t}{g(t+g)},
 \qquad 1 \leq \theta \leq \frac{a}{g},
\end{equation}
\begin{equation}\label{A7}
\frac{1}{\sqrt{x}y} - \frac{1}{\sqrt{t+x}(t+y)} =
 \frac{\theta t}{\sqrt{x}y(t+y)}, \qquad 1 < \theta <
1+\frac{y}{2x},
\end{equation}
\begin{equation}\label{A8}
\frac{1}{\sqrt{xy}z} - \frac{1}{\sqrt{(t+x)(t+y)}(t+z)}
 = \frac{\theta t}{gz\sqrt{(t+x)(t+y)}}, \qquad 1 < \theta <
\frac{a}{g}
 + \frac{g}{z}.
\end{equation}
Finally we have
\begin{equation}\label{A9}
\frac{b}{t^{3/2}(t+b)} < \frac{1}{t^{3/2}} -
\frac{1}{\sqrt{(t+x)(t+y)(t+z)}}
 < \frac{3a}{2t^{3/2}(t+a)},
\end{equation}
where \(a = (x+y+z)/3\) and \(b = \sqrt{3(xy+xz+yz)}/2\), and
\begin{equation}\label{A10}
\frac{t}{g^{3/2}(t+g)} < \frac{1}{\sqrt{xyz}} -
 \frac{1}{\sqrt{(t+x)(t+y)(t+z)}} < \frac{3t}{2g^{3/2}(t+h)},
\end{equation}
where \(g = (xyz)^{1/3}\) and \(3h^{-1} = x^{-1}+y^{-1}+z^{-1}\).

To prove (\ref{A1}) we write
\[
\frac{1}{\sqrt{t}} - \frac{1}{\sqrt{t+x}} = \frac{\sqrt{t+x}-
\sqrt{t}}{\sqrt{t(t+x)}} =
\frac{x}{\sqrt{t(t+x)}(\sqrt{t+x}+\sqrt{t})}
\]
and replace the last denominator factor by either \(2\sqrt{t}\) or
\(2\sqrt{t+x}\).  Interchange of $t$ and $x$ leads from (\ref{A1})
to
(\ref{A2}).  To prove (\ref{A3}) let \(y = \sqrt{1+x/t}\) and write
\[
\frac{t^{3/2}(t+x)}{x}\left(\frac{1}{t^{3/2}}-\frac{1}{(t+x)^{3/2
}}\right)
 = \frac{y^2}{y^2-1}\left(1-\frac{1}{y^3}\right) =
1+\frac{1}{y(y+1)},
\]
which increases from $1$ to $3/2$ as $t$ increases from $0$ to
$\infty$ and
$y$ decreases from $\infty$ to $1$.  Interchange of $t$ and $x$
leads from
(\ref{A3}) to (\ref{A4}).

If the left side of (\ref{A5}) is put over a common denominator, it
suffices to
observe that
\begin{equation}\label{AX}
t+g \leq \sqrt{(t+x)(t+y)} \leq t+a.
\end{equation}
The left inequality is enough to prove the left inequality in
(\ref{A6}). To
prove the right inequality in (\ref{A6}), we define
\[
\phi (t) = (\sqrt{(t+x)(t+y)} - \sqrt{xy})/t
\]
and note that \(\phi(t)\) tends to $a/g$ as \(t \rightarrow 0\) and
to $1$ as
\(t \rightarrow \infty\).  Differentiation shows that $\phi$
decreases
monotonically, because
\[
t^2\sqrt{(t+x)(t+y)}\phi' = -(ta+g^2) + [(ta+g^2)^2 -
t^2(a^2-g^2)]^{1/2}
 \leq 0,
\]
with equality iff \(x=y\).  Because of (\ref{AX}), (\ref{A6})
implies
 (\ref{A6a}).

Equation (\ref{A7}) is proved by solving for $\theta$ and using
(\ref{A2}).
Likewise, (\ref{A8}) is proved by solving for
\[
\theta = \phi(t) + \frac{t}{t+z}
\]
and using the result just established that \(1 \leq \phi(t) \leq
a/g \).

To prove (\ref{A9}) we use Maclaurin's inequality \cite[Thm.
52]{HLP} to find
that
\[
t^3 + 2bt^2 + 4b^2 t/3 < (t+x)(t+y)(t+z) \leq (t+a)^3,
\]
and hence
\begin{equation}\label{AY}
\sqrt{t}( t+b) < \sqrt{(t+x)(t+y)(t+z)} \leq (t+a)^{3/2}.
\end{equation}
Inequality (\ref{A9}) follows from this and (\ref{A3}).

The proof of (\ref{A10}) uses Maclaurin's inequality and the
inequality of
arithmetic and geometric means to get
\[
\frac{t+g}{g} \leq \left[\frac{(t+x)(t+y)(t+z)}{xyz}\right]^{1/3}
=
\left[\left(1+\frac{t}{x}\right)\left(1+\frac{t}{y}\right)
 \left(1+\frac{t}{z}\right)\right]^{1/3}  \leq 1 + \frac{t}{h},
\]
with equalities iff \(x=y=z\), whence
\begin{equation}\label{AZ}
(t+g)^{3/2} \leq \sqrt{(t+x)(t+y)(t+z)} \leq
\left(\frac{g}{h}\right)^{3/2}
 (t+h)^{3/2}.
\end{equation}
Two applications of (\ref{A4}) complete the proof of (\ref{A10}).
\\
\\
\noindent {\bf Acknowledgment.}\quad We thank Arthur Gautesen for
suggesting
the use of uniform approximations.

\end{document}